\ifx\shlhetal\undefinedcontrolsequence\let\shlhetal\relax\fi

\documentclass[11pt]{amsart}

\usepackage{amsmath}
\usepackage{amssymb}

\newtheorem{theorem}{Theorem}[section]
\newtheorem{claim}[theorem]{Claim}

\newtheorem{lemma}[theorem]{Lemma}

\newtheorem{corollary}[theorem]{Corollary}

\theoremstyle{definition}
\newtheorem{definition}[theorem]{Definition}

\newtheorem{question}[theorem]{Question}

\theoremstyle{remark}
\newtheorem{remark}[theorem]{Remark}

\newcount\skewfactor
\def\mathunderaccent#1#2 {\let\theaccent#1\skewfactor#2
\mathpalette\putaccentunder}
\def\putaccentunder#1#2{\oalign{$#1#2$\crcr\hidewidth
\vbox to.2ex{\hbox{$#1\skew\skewfactor\theaccent{}$}\vss}\hidewidth}}

\def\smallbox#1{\leavevmode\thinspace\hbox{\vrule\vtop{\vbox
   {\hrule\kern1pt\hbox{\vphantom{\tt/}\thinspace{\tt#1}\thinspace}}
   \kern1pt\hrule}\vrule}\thinspace}

\newcommand{\cf}{{\rm cf}}

\newcommand{\then}{{\underline{then}}}
\newcommand{\Then}{{\underline{Then}}}

\def\qedref#1{$\qed_{\reforiginal{#1}}$}

\title{The ultrafilter number for singular cardinals}
\author{Shimon Garti}
\address{Institute of Mathematics
 The Hebrew University of Jerusalem
 Jerusalem 91904, Israel}
\email{shimon.garty@mail.huji.ac.il}

\author{Saharon Shelah}
\address{Institute of Mathematics
 The Hebrew University of Jerusalem
 Jerusalem 91904, Israel
 and  Department of Mathematics
 Rutgers University
 New Brunswick, NJ 08854, USA}
\email{shelah@math.huji.ac.il}
\urladdr{http://www.math.rutgers.edu/\char`\~shelah}

\thanks{First typed: December 2011 \newline Research supported by the United States-Israel Binational Science Foundation, and ISF grant no. 242/03. This is publication 1002 of the second author}

\begin{document}
\let\labeloriginal\label
\let\reforiginal\ref

\begin{abstract}
We prove the consistency of a singular cardinal $\lambda$ with small value of the ultrafilter number $\mathfrak{u}_\lambda$, and arbitrarily large value of $2^\lambda$.
\end{abstract}

\maketitle

\newpage

\section{introduction}

Cardinal invariants are defined, traditionally, on the continuum. For some invariants, the continuum is the Cantor space $^\omega 2$, for others it is the Baire space $^\omega \omega$, and sometimes it is the family of infinite sets $[\omega]^\omega$. All these definitions can be generalized to uncountable cardinals, rather than $\omega$. Here we deal with the ultrafilter number, i.e., the minimal cardinality of a generating set for some uniform ultrafilter (the detailed definition is given in the beginning of the next section).

We concentrate on the case of a singular cardinal $\lambda$, aiming to show the consistency of $\mathfrak{u}_\lambda=\lambda^+$ no matter how large is $2^\lambda$. Naturally, for handling successors of singulars we employ methods of $pcf$ theory. We refer to the monograph \cite{MR1318912} and the survey papers \cite{MR2768693},\cite{MR1086455} for general background in this subject. The consistency of the assumptions of the main theorem below is ensured by the forcing construction of \cite{GaSh949}. We quote the pertinent result (this is Claim 3.3 of \cite{GaSh949}):

\begin{theorem}
\label{ffforcing} Product dominating $\lambda$-reals. \newline 
Assume there is a supercompact cardinal in the ground model. \newline 
\Then\ one can force the existence of a singular cardinal $\lambda>\cf(\lambda)=\kappa$, a limit of measurables $\bar{\lambda}=\langle\lambda_i:i<\kappa\rangle$, such that $2^{\lambda_i}=\lambda_i^+$ for every $i<\kappa$ and both products $\prod\limits_{i<\kappa}\lambda_i /J_\kappa^{\rm bd}$ and $\prod\limits_{i<\kappa}\lambda^+_i /J_\kappa^{\rm bd}$ are $\cf(\Upsilon)$-directed for some prescribed $\Upsilon\geq\lambda^{++}$.
\end{theorem}

\hfill \qedref{ffforcing}

\begin{remark}
\label{ccofkappa}
The proof of this theorem in \cite{GaSh949} is phrased for $\lambda>\cf(\lambda)=\omega$, using Prikry forcing. But as indicated in Remark 3.2 there, the same proof can be rendered for the general case of $\lambda>\cf(\lambda)=\kappa>\aleph_0$ using Magidor's forcing. For a survey of Prikry and Magidor forcing, we suggest the excellent work \cite{MR2768695}. So we phrase the theorem in the general setting, but the reader should bare in mind that the full explicit proof for the above theorem is written in \cite{GaSh949} just for singulars with countable cofinality.
\end{remark}

We indicate that the main part in the forcing of \cite{GaSh949} is an iteration of length $\Upsilon$. In the new universe we have ${\rm tcf}(\prod\limits_{i<\kappa}\lambda_i,<_E)={\rm tcf}(\prod\limits_{i<\kappa}\lambda_i^+,<_E)=\cf(\Upsilon)$, so one can choose $\Upsilon$ as a large enough ordinal of cofinality $\lambda^+$ in order to get the desired true cofinality of these products. This important point will be used in the main theorem below.

Our notation is standard. For cardinal invariants we follow \cite{MR2768685}, and the \emph{pcf} notation is due to \cite{MR1318912}. All the filters and ultrafilters in this paper are non-principal. In most cases we also assume that the filter (or ultrafilter) is uniform. It means that each member of the filter has the same cardinality. Suppose $\mathcal{B}$ is a collection of sets. $\mathcal{B}$ has the finite intersection property if $\bigcap\{B_\ell:\ell<n\}$ is not empty for every finite subfamily $\{B_\ell:\ell<n\}$ of $\mathcal{B}$. A collection of sets generates a filter iff it has the finite intersection property. In order to get uniform filters, we add the demand of strong finite intersection property, which means that the size of the intersections is the size of the sets in $\mathcal{B}$.

If $A,B\subseteq\lambda$ then $A\subseteq^*B$ means that $|A\setminus B|<\lambda$.
We denote the diagonal intersection by $\Delta$, so if $\{A_\alpha:\alpha<\kappa\}$ is a collection of subsets of $\kappa$ then $\Delta\{A_\alpha:\alpha<\kappa\}= \{\beta<\kappa: \beta\in\bigcap\limits_{\alpha<\beta}A_\alpha\}$. Recall that if $\lambda$ is a singular cardinal, $\langle\lambda_i:i<\kappa\rangle$ a sequence of regular cardinals which tends to $\lambda$ and $(\prod\limits_{i<\kappa}\lambda_i, <_{J_\kappa^{\rm bd}})$ is $\Upsilon$-directed, then $2^\lambda\geq\Upsilon$. In particular, $2^\lambda\geq\Upsilon$ for arbitrarily large prescribed $\Upsilon$ in Theorem \ref{ffforcing}.

\newpage 

\section{The ultrafilter number}

We commence with the following definition:
\begin{definition}
\label{uuuu} The ultrafilter number. \newline 
Let $\lambda$ be an infinite cardinal, and $\mathcal{F}$ a uniform filter on $\lambda$.
\begin{enumerate}
\item [$(\aleph)$] A base $\mathcal{A}$ for $\mathcal{F}$ is a subfamily of $\mathcal{F}$ such that for every $B\in\mathcal{F}$ there is some $A\in\mathcal{A}$ with the property $A\subseteq^*B$.
\item [$(\beth)$] The ultrafilter number $\mathfrak{u}_\lambda$ is the minimal cardinality of a filter base for some uniform ultrafilter on $\lambda$.
\end{enumerate}
\end{definition}

Let us phrase the following general observation:

\begin{claim}
\label{youlambda} The magnitude of $\mathfrak{u}_\lambda$. \newline 
$\mathfrak{u}_\lambda>\lambda$ for every infinite cardinal $\lambda$.
\end{claim} 

\par\noindent\emph{Proof}.\newline 
Suppose $\mathcal{A}\subseteq\mathcal{P}(\lambda), |\mathcal{A}|\leq\lambda$, $\mathcal{A}$ is closed under finite intersections and endowed with the strong finite intersection property. We try to show that $\mathcal{A}$ does not generate an ultrafilter. For this, we have to designate a subset $Y$ of $\lambda$ so that no member of $\mathcal{A}$ is almost included in $Y$ or in $\lambda\setminus Y$.

We shall define a function $f:\lambda\rightarrow\{0,1\}$ so that:

$$
\forall A\in\mathcal{A}, |A\cap f^{-1}(\{0\})|=|A\cap f^{-1}(\{1\})|=\lambda
$$

If we succeed then we are done. Indeed, set $Y=f^{-1}(\{0\})$ (hence $\lambda\setminus Y=f^{-1}(\{1\})$) and notice that $\neg(A\subseteq^* Y)\wedge\neg(A\subseteq^* \lambda\setminus Y)$ holds for every $Ain\mathcal{A}$. So we need to create the function $f$.

Let $\langle A_i:i<\lambda\rangle$ be an enumeration of the members of $\mathcal{A}$ such that each $A_i$ appears $\lambda$-many times. By induction on $i<\lambda$ we pick an ordinal $\alpha_i$ so that if $i\in\{2j,2j+1\}$ then $\alpha_i\in A_j\setminus\{\alpha_\varepsilon:\varepsilon<i\}$. This can be rendered simply because $|A_j|=\lambda$ while $|\{\alpha_\varepsilon: \varepsilon<i\}|<\lambda$.

Now for every $\alpha<\lambda$ let $f(\alpha)$ be one iff there exists an ordinal $j$ such that $\alpha=\alpha_{2j+1}$, and zero otherwise. Notice that $f(\alpha)$ equals zero whenever $\alpha=\alpha_{2j}$, so by the process of choosing the $\alpha_i$-s we are done.

\hfill \qedref{youlambda}

For one of the assumptions in the main theorem of this paper, we need an ultrafilter which is generated by a $\subseteq^*$-decreasing sequence of sets. We shall use the following:

\begin{lemma}
\label{iincreasing} Almost inclusion generating sequence. \newline 
Suppose $\lambda_i$ is a measurable cardinal, and $U_i$ is a normal ultrafilter on $\lambda_i$. Assume that $2^{\lambda_i}=\lambda_i^+$. 
\then\ there exists a $\subseteq^*$-decreasing sequence $\langle A_{i,\alpha}:\alpha<\lambda_i^+\rangle$ which generates $U_i$.
\end{lemma}

\par\noindent\emph{Proof}.\newline 
Choose an enumeration $\{B_\gamma:\gamma<\lambda_i^+\}$ of $U_i$. For every $\alpha<\lambda_i^+$ let $A_{i,\alpha}$ be $\Delta\{B_\gamma:\gamma<\alpha\}$. By the normality of $U_i$ we know that $A_{i,\alpha}\in U_i$ for every $\alpha<\lambda^+_i$.

By the very definition of the diagonal intersection, if $\alpha_0<\alpha_1<\lambda^+_i$ then $A_{i,\alpha_1}\subseteq^*A_{i,\alpha_0}$, hence the sequence $\langle A_{i,\alpha}:\alpha<\lambda_i^+\rangle$ is $\subseteq^*$-decreasing.
Choose any $B\in U_i$, and let $\gamma$ be an ordinal such that $B\equiv B_\gamma$. For $\alpha=\gamma+1$, $B_\gamma$ appears in the diagonal intersection which defines $A_{i,\alpha}$, so $A_{i,\alpha}\subseteq^*B_\gamma =B$ as required.

\hfill \qedref{iincreasing}

We can state now our main theorem:

\begin{theorem}
\label{mt} Bounding the ultrafilter number. \newline 
Assume that:
\begin{enumerate}
\item [$(\alpha)$] $\kappa=\cf(\lambda)<\lambda$,
\item [$(\beta)$] $E$ is a uniform ultrafilter on $\kappa$,
\item [$(\gamma)$] $\lambda$ is a strong limit cardinal,
\item [$(\delta)$] $\langle\lambda_i:i<\kappa\rangle$ is a sequence of regular cardinals which tends to $\lambda$,
\item [$(\varepsilon)$] $U_i$ is a uniform ultrafilter over $\lambda_i$ for every $i<\kappa$,
\item [$(\zeta)$] For every $i<\kappa$ there is a $\subseteq^*$-decreasing sequence $\langle A_{i,\alpha}:\alpha<\theta_i\rangle$ which generates $U_i$,
\item [$(\eta)$] $\chi_{\bar{\lambda}}={\rm tcf}(\prod\limits_{i<\kappa}\lambda_i,<_E)$ and $\chi_{\bar{\theta}}={\rm tcf}(\prod\limits_{i<\kappa}\theta_i,<_E)$.
\end{enumerate}
\Then\ $\mathfrak{u}_\lambda\leq \chi_{\bar{\lambda}}\cdot\chi_{\bar{\theta}}$.
\end{theorem}

\par\noindent\emph{Proof}. \newline 
For every $i<\kappa$ we fix a sequence $\langle A_{i,\alpha}: \alpha<\lambda_i^+\rangle$ as ensured by Lemma \ref{iincreasing}.
Let $\bar{f}=\langle f_\alpha:\alpha<\chi_{\bar{\lambda}}\rangle$ be a cofinal sequence in $(\prod\limits_{i<\kappa}\lambda_i,<_E)$ and $\bar{g}=\langle g_\beta:\beta<\chi_{\bar{\theta}}\rangle$ a cofinal sequence in $(\prod\limits_{i<\kappa}\theta_i,<_E)$. Denote $\bigcup\{\lambda_j:j<i\}$ by $\lambda_{<i}$. We are trying to define a collection $\mathcal{B}$ of subsets of $\lambda$ so that $|\mathcal{B}|\leq \chi_{\bar{\lambda}}\cdot\chi_{\bar{\theta}}$ and $\mathcal{B}$ generates a uniform ultrafilter.

For each $\alpha<\chi_{\bar{\lambda}},\beta<\chi_{\bar{\theta}}$ and $Y\in E$, set:

$$
B_{\alpha,\beta,Y}=\{\zeta<\lambda: (\exists i\in Y)(\lambda_{<i}\leq\zeta <\lambda_i)\wedge \zeta\in A_{i,g_\beta(i)} \wedge \zeta>f_\alpha(i)\}
$$

Let $\mathcal{B}$ be $\{B_{\alpha,\beta,Y}:\alpha<\chi_{\bar{\lambda}}, \beta<\chi_{\bar{\theta}}, Y\in E\}$. $\mathcal{B}$ is a subset of $\mathcal{P}(\lambda)$ since every $B_{\alpha,\beta,Y}$ is a subset of $\lambda$. The cardinality of $\mathcal{B}$ is as required (i.e., bounded by $\chi_{\bar{\lambda}}\cdot\chi_{\bar{\theta}}$), as $2^\kappa<\chi_{\bar{\lambda}}\cdot\chi_{\bar{\theta}}$.

\emph{Stage A}: $\mathcal{B}$ generates a uniform filter.

For this, it suffices to show that $\mathcal{B}$ has the finite intersection property. Suppose $\{B_{\alpha_\ell,\beta_\ell,Y_\ell}:\ell<n\}$ is a finite sub-collection of $\mathcal{B}$. let $\alpha$ be ${\rm sup}\{\alpha_\ell+1: \ell<n\}$. It means that $f_{\alpha_\ell}<_E f_\alpha$ for every $\ell<n$, hence one can choose a set $Z^1_\ell\in E$ (for every $\ell<n$) so that $i\in Z^1_\ell \Rightarrow f_{\alpha_\ell}(i)<f_\alpha(i)$. Set $Z_1=\bigcap\limits_{\ell<n}Z^1_\ell$, and it follows that $Z_1\in E$.

A similar process can be applied to the `big product' of the $\theta_i$-s. Let $\beta$ be ${\rm sup}\{\beta_\ell+1: \ell<n\}$. By the nature of $\bar{g}$ we have $\ell<n\Rightarrow g_{\beta_\ell}<_E g_\beta$. For each $\ell<n$ choose a set $Z^2_\ell\in E$ such that $i\in Z^2_\ell \Rightarrow g_{\beta_\ell}(i)<g_\beta(i)$. As above, let $Z_2=\bigcap\limits_{\ell<n}Z^2_\ell$, so $Z_2\in E$.

If $i\in Z_2$ then $A_{i,g_\beta(i)}\subseteq^* A_{i,g_{\beta_\ell}(i)}$ for every $\ell<n$, but we need a full inclusion rather than almost inclusion. We can arrange this by eliminating the small exceptions. Fix a pair $(i,\ell)$ so that $i\in Z^2_\ell$, and let $t_{i,\ell}\in [\lambda_i]^{<\lambda_i}$ satisfy $A_{i,g_\beta(i)}\setminus t_{i,\ell} \subseteq A_{i,g_{\beta_\ell}(i)}$. For every $\ell<n$ we define a function $h_\ell\in\prod\limits_{i<\kappa}\lambda_i$ as follows: $h_\ell(i)={\rm sup}(t_{i,\ell})$. Now for every $\ell<n$ choose an ordinal $\gamma_\ell<\chi_{\bar{\lambda}}$ such that $h_\ell<_E f_{\gamma_\ell}$. Define for every $\ell<n$:

$$
Z^3_\ell=\{i<\kappa:A_{i,g_\beta(i)}\setminus f_{\gamma_\ell}(i)\subseteq A_{i,g_{\beta_\ell}(i)}\}
$$

Clearly, $Z^3_\ell\in E$, hence also $Z_3=\bigcap\limits_{\ell<n}Z^3_\ell\in E$. We define $Z=Z_1\cap Z_2\cap Z_3$ and $Y=\bigcap\limits_{\ell<n}Y_\ell\cap Z$, and since we have a finite amount of intersections $Y\in E$. Pick any $i\in Y$. We claim that $\bigcap\{B_{\alpha_\ell,\beta_\ell,Y_\ell}: \ell<n\}\cap[\lambda_{<i},\lambda_i)\in U_i$, in particular this intersection is not empty. Indeed, $B_{\alpha,\beta,Y}\subseteq B_{\alpha_\ell,\beta_\ell,Y_\ell}$ for every $\ell<n$, so $B_{\alpha,\beta,Y}\subseteq \bigcap\{B_{\alpha_\ell,\beta_\ell,Y_\ell}: \ell<n\}$, but $B_{\alpha,\beta,Y}\cap[\lambda_{<i},\lambda_i)\in U_i$ so the finite intersection property is established.

Moreover, since $Y\in E$ we have $\bigcap\{B_{\alpha_\ell,\beta_\ell,Y_\ell}: \ell<n\} \cap[\lambda_{<i},\lambda_i)\in U_i$ for $\kappa$-many $i$-s, and since each $U_i$ is uniform we conclude that the size of this intersection is $\lambda$. It follows that $\mathcal{B}$ generates a uniform filter.

\emph{Stage B}: $\mathcal{B}$ generates an ultrafilter.

Let $X$ be any subset of $\lambda$. For every $i<\kappa$ there is some $t_i\in\{0,1\}$ such that $X\cap\lambda_i\in U_i\Leftrightarrow t_i=1$, so there exists $t\in\{0,1\}$ such that $\{i<\kappa:t_i\equiv t\}\in E$. Without loss of generality $t=1$ (upon replacing $X$ by $\lambda\setminus X$ if unfortunately $t=0$). We denote the set $\{i<\kappa:t_i=1\}$ by $Y_1$.

For every $i\in Y_1$ we choose an ordinal $\beta_i<\theta_i$ so that $A_{i,{\beta_i}}\subseteq^*X\cap\lambda_i$. Clearly, $A_{i,\beta}\subseteq^*X\cap\lambda_i$ for every $\beta\in[\beta_i,\theta_i)$.
We choose an ordinal $\beta<\chi_{\bar{\theta}}$ such that $Y_2=\{i\in Y_1: g_\beta(i)>\beta_i\}\in E$. By the definition of the relation $\subseteq^*$ we know that $\lambda_i>{\rm sup}(A_{i,g_\beta(i)}\setminus X)$ for every $i\in Y_2$. Consequently, one can choose an ordinal $\alpha<\chi_{\bar{\lambda}}$ such that $Y_3=\{i\in Y_2: f_\alpha(i)>{\rm sup}(A_{i,g_\beta(i)}\setminus X)\}\in E$.

Now we can finish the proof by noticing that $B_{\alpha,\beta,Y_3}\subseteq X$. For this, let $\zeta$ be any member of $B_{\alpha,\beta,Y_3}$. It means that $\zeta\in[\lambda_{<i},\lambda_i)$ for some $i\in Y_3$, and $\zeta\in A_{i,g_\beta(i)}\setminus f_\alpha(i)$. By the fact that $i\in Y_3\subseteq Y_1$ we infer that $A_{i,\beta_i}\subseteq^* X\cap\lambda_i$, so $A_{i,\beta}\subseteq^* X\cap\lambda_i$ as well. By the fact that $i\in Y_3$ we have $f_\alpha(i)>{\rm sup}(A_{i,g_\beta(i)}\setminus X)$ so $\zeta\in X$ and the proof is completed.

\hfill \qedref{mt}

We can state now the following corollary:

\begin{corollary}
\label{ssing} Large $2^\lambda$ and small $\mathfrak{u}_\lambda$. \newline 
Suppose there is a supercompact cardinal in ${\rm \bf V}$. \newline 
\then\ it is consistent that there exists a singular cardinal $\lambda$ so that $\mathfrak{u}_\lambda=\lambda^+$ while $2^\lambda$ is arbitrarily large.
\end{corollary}

\par \noindent \emph{Proof}. \newline 
Let $\lambda$ be a singular cardinal, limit of measurables $\langle \lambda_i:i<\kappa\rangle$, $2^{\lambda_i}=\lambda_i^+$ for every $i<\kappa$, and choose any regular cardinal $\tau$ above $\lambda$. Denote $\lambda_i^+$ by $\theta_i$ for every $i<\kappa$. By virtue of Theorem \ref{ffforcing} (and the remarks after it) we may assume that $\chi_{\bar{\lambda}}=\chi_{\bar{\theta}}=\lambda^+$, while $2^\lambda\geq\tau$.

Choose a normal ultrafilter $U_i$ on $\lambda_i$ for every $i<\kappa$. Lemma \ref{iincreasing} provides us with generating $\subseteq^*$-decreasing sequences for every such ultrafilter. It follows that all the requirements in Theorem \ref{mt} are at hand, hence its conclusion holds, namely $\mathfrak{u}_\lambda=\lambda^+$, so the proof is accomplished.

\hfill \qedref{ssing}

Let ${\rm Sp}_\chi(\lambda)$ be the spectrum of cardinals which realize the size of a base for some ultrafilter on $\lambda$. The case of $\lambda=\aleph_0$ was investigated to some extent (see the recent paper \cite{Sh915}). We may wonder what happens in the case of an ubcountable cardinal. In particular, we can ask:

\begin{question}
\label{xxxx} Is it consistent that ${\rm Sp}_\chi(\lambda)$ is not a convex set?
\end{question}

It seems that the methods of proof in this paper might be useful, and we hope to shed light on it in a subsequent work.

\newpage 

\bibliographystyle{amsplain}
\bibliography{arlist}

\end{document}